\documentclass[11pt]{article}
\usepackage{amsmath,amssymb}

\newtheorem{propo}{{\bf Proposition}}[section]
\newtheorem{coro}[propo]{{\bf Corollary}}
\newtheorem{lemma}[propo]{{\bf Lemma}} \newtheorem{theor}[propo]{{\bf
Theorem}} \newtheorem{ex}{{\sc Example}}[section]

\newenvironment{proof}{{\bf Proof.}}{$\Box$}
\def\R{{\mathbb R}}
\def\C{{\mathbb C}}

\begin{document}

\vspace*{1.0in}

\begin{center} SUPPLEMENTS TO MAXIMAL SUBALGEBRAS OF LIE ALGEBRAS 
\end{center}
\bigskip

\begin{center} DAVID A. TOWERS 
\end{center}
\bigskip

\begin{center} Department of Mathematics and Statistics

Lancaster University

Lancaster LA1 4YF

England

d.towers@lancaster.ac.uk 
\end{center}
\bigskip

\begin{abstract}
For a Lie algebra $L$ and a subalgebra $M$ of $L$ we say that a subalgebra $U$ of $L$ is a {\em supplement} to $M$ in $L$ if $L = M + U$. We investigate those Lie algebras all of whose maximal subalgebras have abelian supplements, those that have nilpotent supplements, those that have nil supplements, and those that have supplements with the property that their derived algebra is inside the maximal subalgebra being supplemented. For the algebras over an algebraically closed field of characteristic zero in the last three of these classes we find complete descriptions; for those in the first class partial results are obtained. 
\par
 
\noindent {\em Mathematics Subject Classification 2010}: 17B05, 17B20, 17B30, 17B50.
\par
\noindent {\em Key Words and Phrases}: Lie algebras, maximal subalgebra, supplement, solvable, supersolvable, Frattini ideal. . 
\end{abstract}

\section{Introduction}
\medskip

Let $L$ be a Lie algebra and let $M$ be a subalgebra of $L$. We say that a subalgebra $U$ of $L$ is a {\em supplement} to $M$ in $L$ if $L = M + U$. In similar fashion to \cite{baum} we introduce the following conditions:
\begin{itemize}
\item{(MO)} Every maximal subalgebra of $L$ admits a supplement which is one-dimensional; that is every maximal subalgebra of $L$ has codimension one in $L$.
\item{(MA)} Every maximal subalgebra of $L$ admits an abelian supplement.
\item{(MD)} Every maximal subalgebra $M$ of $L$ admits a supplement whose derived algebra is inside $M$. 
\item{(MN)} Every maximal subalgebra of $L$ admits a nilpotent supplement.
\item{(MU)} Every maximal subalgebra of $L$ admits a supplement every element of which acts nilpotently on $L$. 
\end{itemize}

We will denote by ${\cal MO}$ (respectively, ${\cal MA}, {\cal MD}, {\cal MN}$, and ${\cal MU}$) the class of Lie algebras satisfying condition (MO) (respectively, (MA), (MD), (MN), and (MU)). Our objective is to study these classes of algebras. Corresponding classes of groups were studied by Baumeister (\cite{baum}), who showed, in particular, that any group in ${\cal MA}$ is solvable, and that a group belongs to ${\cal MD}$ if and only if it is solvable. Similar problems concerning factorisations of Lie algebras as sums of subalgebras of a certain type have been studied extensively (see, for example, \cite{o1}, \cite{o2}, \cite{btt}, \cite{gm}, \cite{tow}, \cite{supp}, \cite{cideal}, \cite{index} and the references contained therein.)   
\par

In section two we collect together a few preliminary results. First the description of the algebras in ${\cal MO}$ as derived in \cite{codone} is given. Then some straightforward inclusions between these classes of algebras are noted. Next it is shown that all solvable algebras belong to ${\cal MD}$, that a completely solvable algebra is in ${\cal MU}$ if and only if it is nilpotent, and that all supersolvable and all metabelian algebras are in ${\cal MA}$. A relationship is given between decompositions of a nonassociative algebra and corresponding decompositions of the algebra over a finite field extension. The final result here asserts that if ${\cal H}$ is a saturated homomorph of Lie algebras, then so is $\cal{MH}$.  
\par

The third section is concerned with the simple algebras in these classes. It is shown that if the underlying field is algebraically closed of characteristic zero, then $A_1$ is the only such algebra. We also identify when a maximal parabolic subalgebra has an abelian supplement.
\par

The last section contains the main classification results, describing explicitly the algebras in each of the classes defined above.

\par

Throughout $L$ will denote a finite-dimensional Lie algebra over a field $F$. If $A$ and $B$ are subalgebras of $L$ for which $L = A + B$ and $A \cap B = 0$ we will write $L = A \oplus B$. The ideals $L^{(k)}$ of the derived series are defined inductively by $L^{(0)} = L$, $L^{(k+1)} = [L^{(k)},L^{(k)}]$  for $k \geq 0$; we also write $L^2$ for $L^{(1)}$. We say that $L$ is {\em completely solvable} if $L^2$ is nilpotent. If $A$ is a subalgebra of $L$, the {\em centralizer} of $A$ in $L$ is $C_{L}(A) = \{ x \in L : [x, A] = 0\}$.

\bigskip

\section{Preliminary results}

The Frattini ideal of $L$, $\phi(L)$, is the largest ideal of $L$ contained in
all maximal subalgebras of $L$. The algebras in ${\cal MO}$ were classified in \cite{codone} as follows.

\begin{theor}\label{t:mo} (\cite[Theorem 1]{codone})
Let $L$ be Lie algebra over any field $F$. Then the following are equivalent:
\begin{itemize}
\item[(i)] $L \in {\cal MO}$; and
\item[(ii)] $L/\phi(L) = R \oplus S$ where the radical $R$ is supersolvable and $\phi$-free, and either $S = 0$ or $S$ is three-dimensional simple with a basis $\{u_{-1}, u_0, u_1 \}$ and multiplication $[u_{-1}, u_0] = u_{-1}$, $[u_{-1}, u_1] = u_0$, $[u_0, u_1] = u_1$.
\end{itemize}
\end{theor}

There are some easy relationships between the classes of algebras we have introduced.

\begin{lemma}\label{l:inc} 
\begin{itemize}
\item[(i)] ${\cal MO} \subseteq  {\cal MA} \subseteq {\cal MN}$, ${\cal MA} \subseteq {\cal MD}$ and ${\cal MU} \subseteq {\cal MN}$.
\item[(ii)] If $L$ is solvable, or $F$ has at least $\dim L$ elements, then ${\cal MD} \subseteq {\cal MN}$.
\item[(iii)] If $L$ is completely solvable then ${\cal MN} = {\cal MD}$.
\end{itemize}
\end{lemma}
\begin{proof} (i) These inclusions are straightforward.
\par

\noindent(ii) Suppose that $L \in {\cal MD}$, and let $M$ be any maximal subalgebra of $L$. Then there is a subalgebra $B$ of $L$ such that $L = M + B$ and $B^2 \subseteq M$. Let $C$ be a Cartan subalgebra of $B$ (which, under the stated assumptions, exists) and let $B = C \oplus B_1$ be the Fitting decomposition of $B$ relative to ad\,$C$. Clearly $B_1 \subseteq B^2 \subseteq M$, whence $L = M + C$ and $L \in {\cal MN}$. We have established that ${\cal MD} \subseteq {\cal MN}$.
\par
\noindent(iii) Let $L$ be completely solvable and let $M$ be a maximal subalgebra of $L$ with a nilpotent supplement $U$. Then $U^2 = \phi(U) \subseteq \phi(L) \subseteq M$, by \cite{stit}. It follows that ${\cal MN} \subseteq {\cal MD}$. The reverse inclusion comes from (ii) above.
\end{proof}
\bigskip

We define the {\em abelian socle} of $L$, $Asoc\,L$, to be the sum of the minimal abelian ideals of $L$. Next we consider the solvable algebras in these classes.

\begin{propo}\label{p:solv} Let $L$ be a Lie algebra over an arbitrary field $F$. 
\begin{itemize}
\item[(i)] If $L$ is solvable then $L \in {\cal MD}$.
\item[(ii)] If $L$ is completely solvable, then $L \in {\cal MU}$ if and only if $L$ is nilpotent.
\item[(iii)] If $L$ is supersolvable then $L \in {\cal MA}$.
\item[(iv)] If $L$ is metabelian (so, in particular, if $L$ is completely solvable and $\phi$-free) then $L \in {\cal MA}$.
\end{itemize}
\end{propo}
\begin{proof} (i) Let $L$ be solvable and let $M$ be a maximal subalgebra of $L$. Then there is a $k \geq 0$ such that $L^{(k)} \not \subseteq M$ but $L^{(k+1)} \subseteq M$. Clearly $L = M + L^{(k)}$ and so $L \in {\cal MD}$.
\par

\noindent (ii) Let $L$ be completely solvable. Suppose that $L \in {\cal MU}$, but that $L$ is not nilpotent, and let $M$ be a maximal subalgebra of $L$ with $N \subseteq M$, where $N$ is the nilradical of $L$. Then there is a nil subalgebra $U$ of $L$ such that $L = M + U$. But now $N + U$ is nilpotent, and is an ideal of $L$, since $L^2 \subseteq N$, so $U \subseteq N \subseteq M$, a contradiction. The converse is clear. 
\par

\noindent (iii) If $L$ is supersolvable then $L \in {\cal MO} \subseteq {\cal MA}$, by Theorem \ref{t:mo} and Lemma \ref{l:inc}(i).
\par

\noindent (iv) Let $L$ be metabelian and suppose that $M$ is a maximal subalgebra of $L$. If $L^2 \subseteq M$ then $M$ has codimension one in $L$ and so has an abelian supplement. If $L^2 \not \subseteq M$, then $L^2$ is an abelian supplement to $M$. That completely solvable $\phi$-free algebras lie in this class follows from \cite[Theorems 7.3 and 7.4]{frat}. 
\end{proof}
\bigskip

Not every completely solvable Lie algebra $L$ belongs to ${\cal MA}$, as the next example shows.

\begin{ex}\label{e:ex} Let $L$ be the four-dimensional Lie algebra over the real field with basis $e_1, e_2, e_3, e_4$ and products $[e_1,e_2] = e_2 - e_3$, $[e_1, e_3] = e_2 + e_3$, $[e_1,e_4] = 2e_4$, $[e_2,e_3] = e_4$, other products being zero.
\end{ex}
Then $L$ is completely solvable, but not $\phi$-free, as $\phi(L) = \R e_4$. Also, $M = \R e_1 + \R e_4$ is a maximal subalgebra of $L$ that has no abelian supplement.
\bigskip

We have the following relationship between decompositions of a nonassociative algebra and corresponding decompositions of the algebra over a finite field extension.

\begin{lemma}\label{l:fe} Let $A$ be a nonassociative algebra over a field $F$ with subalgebras $A_1$ and $A_2$, and let $K$ be a finite field extension of $F$ so that $F$ is the fixed field of the group $Gal(K/F)$ of $F$-automorphisms of $K$. Then $A = A_1 + A_2$ if and only if $\bar{A} = \bar{A_1} + \bar{A_2}$, where $\bar{B} = B \otimes_F K$ for any subalgebra $B$ of $A$.
\end{lemma}
\begin{proof} Clearly $A = A_1 + A_2$ implies that $\bar{A} = \bar{A_1} + \bar{A_2}$. Conversely, let $\bar{A} = \bar{A_1} + \bar{A_2}$, and let $x \in A$. Then $x \otimes 1 = \bar{x_1} + \bar{x_2}$ where $\bar{x_1} \in \bar{A_1}$, $\bar{x_2} \in \bar{A_2}$. Let $\{k_1, \ldots , k_n\}$ be a basis for $K$ over $F$, $Gal(K/F) = \{\theta_1, \ldots , \theta_n\}$. For each $\theta \in Gal(K/F)$ let $U_{\theta}$ be the semilinear transformation of $\bar{A}$ defined by $U_{\theta}(\sum a_i \otimes k_i) = \sum a_i \otimes \theta(k_i)$ as in \cite[page 295]{jac}. Then $U_{\theta_i}(x \otimes 1) = x \otimes 1$, and so $x \otimes 1 = U_{\theta_i}(\bar{x_1}) + U_{\theta_i}(\bar{x_2})$ for all $1 \leq i \leq n$. Hence  
\[
x \otimes 1 = \frac{1}{n}(U_{\theta_1}(\bar{x_1}) + \ldots + U_{\theta_n}(\bar{x_1})) + \frac{1}{n}(U_{\theta_1}(\bar{x_2}) + \ldots + U_{\theta_n}(\bar{x_2})).
\]
Moreover, $\frac{1}{n}(U_{\theta_1}(\bar{x_i}) + \ldots + U_{\theta_n}(\bar{x_i}))$ is fixed by all of the elements of $\{U_{\theta}:\theta \in Gal(K/F)\}$, and so belongs to $A_i$ for each $i = 1,2$, whence the result.
\end{proof}
\bigskip

Notice, however, that if $L$ is as in Example \ref{e:ex} above then $L \notin {\cal MA}$, whereas, considered as an algebra over $\C$, we have $L \in {\cal MA}$.
\par 

A class ${\cal H}$ of finite-dimensional Lie algebras is called a {\em homomorph} if it contains, along with an algebra $L$, all epimorphic images of $L$; it is {\em saturated} if $L/\phi(L) \in {\cal H}$ implies that $L \in {\cal H}$. Then we have that if ${\cal H}$ is a saturated homomorph, so is ${\cal MH}$ (where ${\cal MH}$ is the class of Lie algebras all of whose maximal subalgebras have a supplement $U \in {\cal H}$.) First we need a lemma.

\begin{lemma}\label{l:min} Let $L$ be a Lie algebra over any field $F$ and let $M$ be a maximal subalgebra of $L$ with a supplement $U$. Then $M$ has a supplement $W$ with $\phi(L) \cap W \subseteq \phi(W)$.
\end{lemma}
\begin{proof} Choose $W$ to be a subalgebra of $U$ that is minimal with respect $U = \phi(L) \cap U + W$. Then $L = M + U = M + \phi(L) \cap U + W = M + W$ and $\phi(L) \cap W = (\phi(L) \cap U) \cap W \subseteq \phi(W)$, by \cite[Lemma 7.1]{frat}. 
\end{proof}
\bigskip  

\begin{propo}\label{p:sh} Let ${\cal H}$ be a saturated homomorph of Lie algebras. Then $L \in \cal{MH}$ if and only if $L/\phi(L) \in \cal{MH}$; that is $\cal{MH}$ is also a saturated homomorph.  
\end{propo}
\begin{proof} Suppose first that $L \in \cal{MH}$ and let $M/\phi(L)$ be maximal subalgebra of $L/\phi(L)$. Then $M$ is a maximal subalgebra of $L$ and so there is a subalgebra $U \in {\mathcal H}$ such that $L = M + U$. But now 
$$\frac{L}{\phi(L)} = \frac{M}{\phi(L)} + \frac{U + \phi(L)}{\phi(L)} \hbox{ and } \frac{U + \phi(L)}{\phi(L)} \cong \frac{U}{U \cap \phi(L)} \in {\mathcal H},$$
whence $L/\phi(L) \in {\cal MH}$.
\par

 So suppose now that $L/\phi(L) \in \cal{MH}$ and let $M$ be a maximal subalgebra of $L$. Then there is a subalgebra $U/\phi(L) \in {\cal H}$ such that $L/\phi(L) = M/\phi(L) + U/\phi(L)$, and so $L = M + U$. Thus there is a subalgebra $W$ of $U$ with $U = \phi(L) + W$, $L = M + W$ and $\phi(L) \cap W \subseteq \phi(W)$, by Lemma \ref{l:min}. Moreover,
$$\frac{W}{\phi(L) \cap W} \cong \frac{\phi(L) + W}{\phi(L)} = \frac{U}{\phi(L)} \in {\mathcal H}, \hbox{ so } \frac{W}{\phi(W)} \cong \frac{W/(\phi(L) \cap W)}{\phi(W)/(\phi(L) \cap W)} \in {\mathcal H}.$$ 
It follows that $W \in {\cal H}$ and therefore $L \in {\cal MH}$.
\par

It is easy to see that ${\cal MH}$ is a homomorph.    
\end{proof}

\section{Simple algebras}
\medskip
Our first objective in this section is to establish the following result.

\begin{theor} The only simple Lie algebra over an algebraically closed field of characteristic zero belonging to ${\cal MA}$, or to ${\cal MU}$, is $A_1$.
\end{theor}
\bigskip

First let us recall some facts about maximal subalgebras of simple Lie algebras $L$ over an algebraically closed fields $F$ of characteristic zero. They fall into three types: 
\begin{itemize}
\item[(I)] reducible maximal subalgebras, which are described in \cite[Theorems 1.1, 1.2, page 252]{dynk2};
\item[(II)] irreducible non-simple maximal subalgebras, which are described in \cite[Theorems 1.3, 1.4, page 253]{dynk2}; and
\item[(III)] irreducible simple maximal subalgebras, which are described in \cite[Theorem 1.5, page 252]{dynk2}.
\end{itemize}
A subalgebra $B$ of a semisimple Lie algebra $L$ is called {\em regular} in $L$ if we can choose a basis for $B$ in such a way that any vector of this basis is either a root vector of $L$ corresponding to some Cartan subalgebra $C$ of $L$, or otherwise belongs to $C$; $B$ is an {\em $R$-subalgebra} of $L$ if it is contained in a regular subalgebra of $L$, and is an {\em $S$-subalgebra} otherwise (see \cite[page 158]{dynk1}). We say that $B$ is {\em parabolic} in $L$ if it contains a Borel subalgebra of $L$; it is {\em reductive} in $L$ if the representation $x \mapsto \hbox{ad}_L\,x$ of $B$ is semisimple.
\par

Then a further way of describing the maximal subalgebras of $L$ is that they are either 
\begin{itemize}
\item[(a)] parabolic, all of which are regular and of type (I); or
\item[(b)] reductive, which further subdivide into:
\begin{itemize}
\item[(i)] regular reductive subalgebras, which are of type (I), and are semisimple of maximal rank; and
\item[(ii)] $S$-subalgebras, which are semisimple and are mostly of type (II) or type (III) (entirely so in the case of $A_n$, $B_n$ and $C_n$).
\end{itemize}
\end{itemize}

\begin{table}
\begin{center}
\begin{tabular}{|c|c||c|c|c|c|} \hline
simple  & maximal  & dim\,$L$ & dim\,$M$ & rank\,$L$ & dim\,$A$ \\ 
algebra  & subalgebra  & & &  &\\ 
 $L$ & $M$  &   &    &    &  dim\,$U$ \\ \hline \hline
$A_2$ & $A_1$ & $8$ & $3$ & $2$ & $2$ \\  
  &   &   &    &    &  $3$ \\\hline 
$A_{2n}$  & $B_n$ & $4n^2 + 4n$ & $2n^2 + n$ & $2n$ & $n^2 + n$ \\
 ($n \geq 2$) &   &   &    &    &  $2n^2 + n$ \\  \hline
$A_{2n+1}$ & $D_{n+1}$ & $4n^2 + 8n + 3$ & $2n^2 + 3n + 1$ & $2n + 1$ & $n^2 + 2n + 1$ \\  
 ($n \geq 1$) &   &   &    &    &  $2n^2 + 3n + 1$ \\\hline 
$B_2$ & $A_1 \oplus A_1$ & $10$ & $6$ & $2$ & $3$ \\
  & $A_1^{10}$  &   &  $3$  &    &  $4$ \\ \hline
$B_3$ & $A_1 \oplus A_1 \oplus A_1$ & $21$ & $9$ & $3$ & $5$ \\
  &   &   &    &    &  $9$ \\  \hline 
$B_{2n}$  & $B_n \oplus D_n$ & $8n^2 + 2n$ & $4n^2$ & $2n$ & $2n^2 - n + 1$ \\
 ($n \geq 2$) &   &   &    &    &  $4n^2$ \\  \hline
$B_{2n+1}$  & $B_n \oplus D_{n+1}$ & $8n^2 + 10n + 3$ & $4n^2 + 4n + 1$ & $2n + 1$ & $2n^2 + n + 1$ \\  
 ($n \geq 2$) &   &   &    &    &  $4n^2 + 4n + 1$ \\\hline 
$C_{2n}$  & $C_n \oplus C_n$ & $8n^2 + 2n$ & $4n^2 + 2n$ & $2n$ & $2n^2 + n$ \\
 ($n \geq 2$) & $A_1$  &   & $3$   &    &  $4n^2$ \\  \hline
$C_{2n+1}$  & $C_n \oplus C_{n+1}$ & $8n^2 + 10n + 3$ & $4n^2 + 6n + 3$ & $2n + 1$ & $2n^2 + 3n + 1$ \\  
 ($n \geq 1$) & $A_1$  &   & $3$   &    &  $4n^2 + 4n + 1$ \\\hline 
$D_{2n}$  & $D_n \oplus D_n$ & $8n^2 - 2n$ & $4n^2 - 2n$ & $2n$ & $2n^2 - n$ \\  
 ($n \geq 2$) &   &   &    &    &  $4n^2 - 2n$ \\\hline
$D_{2n+1}$  & $D_n \oplus D_{n+1}$ & $8n^2 + 6n + 3$ & $4n^2 + 2n + 1$ & $2n + 1$ & $2n^2 + n$ \\  
 ($n \geq 2$) &   &   &    &    &  $4n^2 + 2n + 1$ \\\hline
$E_6$ & $A_2^{9}$ & $78$ & $8$ & $6$ & $16$ \\
  &   &   &    &    &  $36$ \\  \hline 
$E_7$ & $A_1^{231}, A_1^{399}$ & $133$ & $3$ & $7$ & $27$ \\
  &   &   &    &    &  $63$ \\  \hline 
$E_8$ & $A_1^{520}, A_1^{760}$,  & $248$ & $3$ & $8$ & $36$ \\
  & $A_1^{1240}$  &   &    &    &  $120$ \\  \hline 
$F_4$ & $A_1^{156}$ & $52$ & $3$ & $4$ & $9$ \\  
  &   &   &    &    &  $24$ \\\hline 
$G_2$ & $A_1^{28}$ & $14$ & $3$ & $2$ & $3$ \\
  &   &   &    &    &  $6$ \\  \hline 
\end{tabular}
\end{center}
\caption{Maximal Subalgebras}
\end{table}
\bigskip

\noindent{\bf Proof of Theorem 3.1.} For each class of simple Lie algebras $L$ we shall show that there is a maximal subalgebra $M$ such that for any abelian subalgebra $A$ (respectively, nil subalgebra $U$) of maximal dimension in $L$, dim\,$(M + A) <$ dim\,$L$ (respectively, dim\,$(M + U) <$ dim\,$L$). This is the content of the Table 1 below, whose entries we will explain next. 
\par

For each class of algebra we list in column 2 a possible choice of maximal subalgebra $M$ to meet our claim. In some cases two subalgebras are given: the top one is sufficient to rule out the possibility of an abelian supplement, but not of a nil supplement. Of course the lower one of the two would suffice on its own, but the top one is a more straightforward example and so is listed for interest.
\par

Many maximal subalgebras can be found from Dynkin's trick of removing a node from the extended Dynkin diagram: this gives regular reductive subalgebras. The top entry for the algebras of types $B, C$ or $D$ can be found in this way (and remembering that $B_1 = A_1$ and $D_2 = A_1 \oplus A_1$); they can also be found in \cite[Table 12, page 150; page 232]{dynk1}. The lower entries for $C_n$ and $B_2$ follow from the fact that all of the three-dimensional $S$-subalgebras of $B_n$ and $C_n$ are maximal, except for $A_1^{28} \subset G_2 \subset B_3$ (whereas all three-dimensional $S$-subalgebras of $A_n$ and $D_n$ are non-maximal except for $A_1^{4} \subset A_2$.) The superfixes indicate the index of the embedding as defined in \cite{dynk1}.
\par

The three-dimensional representation of $A_1$ gives an embedding of $A_1$ in $A_2$ under which it is a maximal $S$-subalgebra of that algebra. There are natural embeddings of $B_n$ in $A_{2n} (n \geq 2)$ and $D_{n+1}$ in $A_{2n+1} (n \geq 1)$ (also in \cite[Table 5, page 366]{dynk2}) under which these are maximal $S$-subalgebras of those algebras. Finally the maximal $S$-subalgebras listed for each of the exceptional simple Lie algebras are taken from \cite[Table 39, page 233]{dynk1}. 
\par

In the final column of the table the upper number, $\alpha$ (respectively, lower number $\gamma$), is the maximal possible dimension of an abelian (respectively, nil subalgebra). The values for $\alpha$ were determined by Malcev in \cite{mal}, or are given in \cite[Table 1]{bc}. The value for $\gamma$ is computed as $\frac{1}{2}(\hbox{dim}\,L - \hbox{rank}\,L)$ (see \cite{dkk}).
\bigskip

The above result can be extended to cover ${\cal MN}$ as follows.

\begin{coro}\label{c:mn} The only simple Lie algebra over an algebraically closed field of characteristic zero belonging to ${\cal MN}$ is $A_1$. 
\end{coro} 
\begin{proof} Suppose that $L \not \cong A_1$ and let $M$ be a maximal subalgebra of $L$ from Table 1 (the lower entry if two are given) and suppose that $L = M + N$ where $N$ is a maximal nilpotent subalgebra of $L$. Then $N = C_L(N)$ and so $N$ is algebraic. It follows that $N = T \oplus U$ where $T$ is a toral subalgebra and $U$ is the nilradical of $N$. Since $N$ is nilpotent we must have that $[T,U] = 0$. 
\par

Now it can be seen from Table 1 that $N$ must have the same dimension as a Borel subalgebra of $L$. In fact, in view of \cite{dkk}, it must be equal to a Borel subalgebra $B$ of $L$ in which $T$ is a Cartan subalgebra of $L$ and $U$ is the nilradical of $B$. Since $[T,U] = 0$ this is impossible.   
\end{proof}
\bigskip

All of the subalgebras in Table 1, of course, are reductive. The maximal parabolic subalgebras generally are too large to yield to dimension arguments of the above kind. However, we can determine which of them have abelian supplements.

Let $L$ be a simple Lie algebra over an algebraically closed field of characteristic zero, $H$ a Cartan subalgebra of $L$ and $\Gamma$, $\Gamma^+$, $\Sigma$, respectively, the systems of roots, positive roots and simple roots of $L$ with respect to $H$. Then $L = H + \Sigma_{\alpha \in \Gamma} V_{\alpha}$, where $V_{\alpha}$ is spanned by a unique element $e_{\alpha}$. Let $\Sigma_1 \subseteq \Sigma$ be a non-empty subsystem of $\Sigma$, and put $\Delta_1 = \{\gamma \in \Gamma : \gamma = \Sigma_{\alpha \in \Sigma \setminus \Sigma_1} m_{\alpha} \alpha\}$ and $\Delta_2^+ = (\Gamma \setminus \Delta_1) \cap \Gamma^+$. 
\par

Then every parabolic subalgebra of $L$ is conjugate to a standard parabolic subalgebra of the form $P = H + \Sigma_{\alpha \in \Delta_1} V_{\alpha} + \Sigma_{\alpha \in \Delta_2^+} V_{\alpha} = R \oplus U$, where $R = H + \Sigma_{\alpha \in \Delta_1} V_{\alpha}$ is its reductive summand and the ideal $N = \Sigma_{\alpha \in \Delta_2^+} V_{\alpha}$ is its nilradical. It is clear that every maximal parabolic subalgebra has a nil supplement, namely the nilradical, $N^o = \Sigma_{- \alpha \in \Delta_2^+} V_{\alpha}$, of the opposite parabolic subalgebra of $L$. When they have an abelian supplement is given by the next result, where we use the Bourbaki numbering of roots (see \cite{bou}).

\begin{propo} Let $L$ be a simple Lie algebra over an algebraically closed field $F$ of characteristic zero, and let $P$ be a standard maximal parabolic subalgebra of $L$. Then $P$ has an abelian supplement in $L$ if and only if $L \cong A_n$ and $\Sigma_1 = \{\alpha_i\}$, $B_n$ and $\Sigma_1 = \{\alpha_1\}, C_n$ and $\Sigma_1 = \{\alpha_n\}$, $D_n$ and $\Sigma_1 = \{\alpha_1\}, \{\alpha_{n-1}\}$ or $\{\alpha_n\}$, $E_6$ and $\Sigma_1 = \{\alpha_1\}, \{\alpha_6\}$, or $E_7$ and $\Sigma_1 = \{\alpha_7\}$.
\end{propo}
\begin{proof} We have that $L = P \oplus N^o$ where $N_o$ is the nilradical of the opposite parabolic subalgebra to $P$. Suppose that $N^o$ is not abelian, but that $L = P + A$, where $A$ is abelian. We can embed $A$ in a Borel subalgebra $B$ of $L$. Then $B$ is conjugate to a Borel subalgebra $B^o$ containing $N^o$ as its nilradical. Let $A^o$ be the image of $A$ under the conjugating automorphism. Let $\alpha(B^o)$, respectively $\beta(B^o)$, be the maximal possible dimension of an abelian subalgebra, respectively ideal, of $B^o$. Then $\alpha(B^o) = \beta(B^o)$, by \cite[Proposition 2.5]{bc}, so dim\,$A^o \leq \alpha(B^o) = \beta(B^o) <$ dim\,$N^o$ if $N^o$ is not abelian. It follows that $P$ has an abelian supplement precisely when $N^o$ is abelian. But this occurs exactly in the situations given in the result (see \cite{oh}, or \cite[Table 1, page 24]{bajo}).
\end{proof}

\bigskip

\section{Main results}

First we need to extend the results of the previous section to semisimple Lie algebras.

\begin{lemma}\label{l:ss} Let $L$ be a semisimple Lie algebra over an algebraically closed field $F$ of characteristic zero. Then $L \in {\cal MN} = {\cal MA} = {\cal MD}$ if and only if $L \cong A_1$.
\end{lemma}
\begin{proof} Suppose that $L \in {\cal MN}$ is semisimple. Clearly every simple summand of $L$ is isomorphic to $A_1$, by Corollary \ref{c:mn}. Suppose that $L = S \oplus \bar{S}$, where $S \cong A_1$, $\bar{S}$ is an isomorphic copy of $S$ with $[S, \bar{S}] = 0$, and denote the image of $s \in S$ in $\bar{S}$ by $\bar{s}$. Let $D = \{s + \bar{s} : s \in S\}$, which is easily seen to be a maximal subalgebra of $L$, and suppose that $L = D + N$, where $N$ is a nilpotent subalgebra of $L$. Then dim\,$N \geq 3$.
\par

Clearly $L \neq S + N$, so $S \cap N \neq \{0\}$. Similarly, $\bar{S} \cap N \neq \{0\}$.  Let $s \in S \cap N$, $\bar{x} \in \bar{S} \cap N$ and let $n = u + \bar{v} \in N$, where $u \in S$, $\bar{v} \in \bar{S}$. Then $[s, u] = [s, n] \in S \cap N$, so $[s, n] = \lambda s$ for some $\lambda \in F$, since $S$ has no two-dimensional nilpotent subalgebras. As $N$ is nilpotent, $\lambda = 0$. But now $Fs + Fu$ is an abelian subalgebra of $S$, and so $u = \mu s$ for some $\mu \in F$. In similar manner $\bar{v} = \nu \bar{x}$ for some $\nu \in F$. But this means that dim\,$N \leq 2$, a contradiction. It follows that $L$ is simple and $L \cong A_1$. The same conclusion follows if $L \in {\cal MA}$ or $L \in {\cal MD}$, by Lemma \ref{l:inc}.
\par

The converse is easily checked.   
\end{proof}

\begin{theor}\label{t:md} Let $L$ be a Lie algebra over an algebraically closed field of characteristic zero with solvable radical $R$. Then $L \in {\cal MD} = {\cal MN}$ if and only if $L$ is solvable or $L/R \cong A_1$. 
\end{theor}
\begin{proof} Suppose first that $L \in {\cal MN}$ and $L$ is not solvable. Then $L/R \cong A_1$ by Lemma \ref{l:ss}. If $L \in {\cal MD}$ the same conclusion follows from Lemma \ref{l:inc}.
\par

Suppose now that $L = R$ or $L/R \cong A_1$, and let $M$ be a maximal subalgebra of $L$. If $R \subseteq M$ then $M$ has a supplement whose derived algebra is inside $M$, by Lemma \ref{l:ss}. So suppose that $R \not \subseteq M$. Then there is a $k \geq 0$ such that $R^{(k)} \not \subseteq M$ but $R^{(k+1)} \subseteq M$. Clearly $L = M + R^{(k)}$ and again $M$ has a supplement whose derived subalgebra is inside $M$. It follows that $L \in {\cal MD}$. Lemma \ref{l:inc}(ii) also implies that $L \in {\cal MN}$. 
\end{proof}
\bigskip

Notice that Lemma \ref{l:ss} and Proposition \ref{p:solv} imply that if $L$ is a semisimple or solvable Lie algebra over an algebraically closed field of characteristic zero then $L \in {\cal MA}$ if and only if $L \in {\cal MN}$. However, it is not the case that ${\cal MA} = {\cal MN}$, as the following example shows.
\bigskip

\begin{ex} Let $L$ be the six-dimensional Lie algebra over the complex field with basis $e, f, h, x_0, x_1, x_2$ and products $[e, h] = 2e$, $[f, h] = -2f$, $[e, f] = h$, $[x_0, h] = x_0$, $[x_1, h] = -x_1$, $[x_0, f] = x_1$, $[x_1, e] = -x_0$, $[x_0, x_1] = x_2$, other products being zero. 
\end{ex}
Clearly, $L = R \oplus S$, where $R = \C x_0 + \C x_1 + \C x_2$ is nilpotent and $S = \C e + \C f + \C h \cong A_1$. Then $M = \C e + \C f + \C h + \C x_2$ is a maximal subalgebra of $L$ that has no abelian supplement in $L$.

\begin{theor} Let $L$ be a Lie algebra over an algebraically closed field of characteristic zero with nilradical $N$. Then $L \in {\cal MU}$ if and only if $L$ is nilpotent or $L/N \cong A_1$. 
\end{theor}
\begin{proof} Let $L \in {\mathcal MU}$. Then $L \in {\mathcal MN}$, by Lemma \ref{l:inc}, and so $L$ is solvable or $L/R \cong A_1$, by Theorem \ref{t:md}. Suppose that $L$ is solvable but not nilpotent, and let $M$ be a maximal subalgebra of $L$ with $N \subseteq M$. Then there is a nil subalgebra $U$ of $L$ such that $L = M + U$. But now $N + U$ is nilpotent, and is an ideal of $L$, since $L^2 \subseteq N$, so $U \subseteq N \subseteq M$, a contradiction. 
\par

So suppose now that $L/R \cong A_1$ and $R$ is not nilpotent. Clearly $L/\phi(L) \in {\mathcal MU}$, so assume that $L$ is $\phi$-free. Then $L = Asoc\,L \oplus (C \oplus S)$, where $C$ is abelian and acts semisimply on $Asoc\,L$, $S \cong A_1$, and $[S,C] = 0$, by \cite[Theorem 7.5]{frat}. Let $M$ be a maximal subalgebra of $L$ with $Asoc\,L + S \subseteq M$. There is a nil subalgebra $U$ such that $L = M + U$. Since $L^2 \subseteq M$, $M$ is an ideal of $L$ and there is a $u \in U$ such that $L = M + Fu$. Let $u = a + c + s$, where $a \in Asoc\,L$, $c \in C$, $s \in S$. It is easy to see that since ad\,$(a+c+s)$ acts nilpotently on $S$, then $s$ must be a nil element of $S$. But now ad\,$(a+c+s)|_{Asoc\,L} =$ ad\,$(c+s)|_{Asoc\,L}$ is nilpotent. As $c$ acts semisimply on $Asoc\,L$ it follows that $c \in C_L(Asoc\,L) \subseteq Asoc\,L$. But then $u \in M$, a contradiction. Hence $R$ is nilpotent. 
\par

Suppose conversely that $L$ is nilpotent or $L/N \cong A_1$. If the former holds that then, clearly, $L \in {\mathcal MU}$. So suppose the latter holds and let $M$ be a maximal subalgebra of $L$. If $N \not \subseteq M$ then $L = M + N$ and $N$ is nil. If $N \subseteq M$ then $M = N + M \cap A_1$. Also $M \cap A_1$ is a maximal subalgebra of $A_1$ and there is a nil element of $A_1$, $e$ say, such that $L = M + Fe$. But then $e$ acts nilpotently on $L$. It follows that $L \in {\mathcal MU}$.
\end{proof}


\bigskip


\begin{thebibliography}{1}

\bibitem{baum} {\sc B. Baumeister}, `A characterisation of finite soluble groups', {\em Arch. Math.} {\bf 72} (1999), 167--176.

\bibitem{btt} {\sc Y. Bahturin, M. Tvalavadze and T. Tvalavadze}, `Sums of simple and nilpotent Lie algebras', {\em Comm. Alg.} {\bf 30} (2002), 4455--4471.

\bibitem{bajo} {\sc I. Bajo}, `Maximal isotropy groups of Lie groups related to nilradicals of parabolic subalgebras', {\em Tohoku Math. J.} {\bf 52} (2000), 19--29.

\bibitem{bou} {\sc N. Bourbaki}, `Groupes et algebres de Lie, Chapitres 4, 5 et 6', Hermann, Paris 1968, Masson, Paris 1981; `Lie groupes and Lie algebras, Chapters 4-6', Translated from the 1968 French original by Andrew Pressley, Springer, Berlin 2002.

\bibitem{bc} {\sc D. Burde and M. Ceballos}, `Abelian ideals of maximal dimension for solvable Lie algebras', {\em ArXiv: 0911.2995v1 [math.RA]} (2009).

\bibitem{dkk} {\sc J. Draisma, H. Kraft and J. Kuttler}, `Nilpotent subspaces of maximal dimension in semisimple Lie algebras', {\em Comp. Math.} {\bf 142} (2006), 464--476.

\bibitem{dynk1} {\sc E. B. Dynkin}, `Semisimple subalgebras of semisimple Lie algebras', {\em Mat. Sbornik} {\bf 30}(72), (1952), 349 -- 462; {\em Amer. Math. Soc. Transl. (2)} {\bf 6} (1957), 111--244. 

\bibitem{dynk2} {\sc E. B. Dynkin}, `The maximal subgroups of the classical groups', {\em Proc. Moscow Math. Soc.} {\bf 1}, (1952), 39 -- 166; {\em Amer. Math. Soc. Transl. (2)} {\bf 6} (1957), 245--378. 

\bibitem{gm} {\sc A.G. Gein and Yu.N. Mukhin}, `Complements to subalgebras of Lie algebras', {\em Ural. Gos. Univ. Mat. Zap.} {\bf 12}, No. 2 (1980) 24--48.

\bibitem{jac} Jacobson, N.: Lie algebras. New York: Dover Publ. (1979).

\bibitem{mal} {\sc A. Malcev}, `Commutative subalgebras of semisimple Lie algebras', {\em Bull. Acad. Sci. URSS Ser. Math.} {\bf 9} (1945), 291 -- 300 [Izvestia Akad. Nauk SSSR]; {\em Amer. Math. Soc. Transl.} {\bf 40} (1951).

\bibitem{o1} {\sc A.L. Onishchik}, `Inclusion relations among transitive compact transformation groups', {\em Trudy Moskov. Mat. Obshch.} {\bf 11} (1962), 199--242; {\em Amer. Math. Soc. Transl.}(2) {\bf 50} (1966), 5--58.

\bibitem{o2} {\sc A.L. Onishchik}, `Decompositions of reductive Lie groups', {\em Mat. Sbornik} {\bf 80}(122)(4) (1969), 515--554; {\em Math. USSR-Sbornik} {\bf 9}(4) (1969), 515--554. 

\bibitem{oh} {\sc A.L. Onishchik and Yu.B. Hakimjanov}, `Semidirect sums of Lie algebras', {\em Mat. Zametki} {\bf 18} (1975), 31--40; {\em Math. Notes} {\bf 18} (1975), 600--604. 

\bibitem{stit} {\sc E.L. Stitzinger}, `Frattini subalgebras of a class of solvable Lie algebras', {\em Pacific J. Math.} {\bf 34}, No. 1 (1970), 177--182.

\bibitem{frat} {\sc D.A. Towers}, `A Frattini theory for algebras', {\em
Proc. London Math. Soc.} (3) {\bf 27} (1973), 440--462.

\bibitem{tow} {\sc D.A. Towers}, `On Lie algebras in which modular pairs of subalgebras are permutable', {\em J. Algebra} {\bf 68}, No. 2 (1981), 369--377.

\bibitem{codone} {\sc D.A. Towers}, `Lie algebras, all of whose maximal subalgebras have codimension one', {\em Proc. Edin. Math. Soc.} {\bf 24} (1981), 217--219.

\bibitem{supp} {\sc D.A. Towers}, `C-supplemented subalgebras of Lie algebras', {\em J. Lie Theory} {\bf 18}, (2008), 717--724.

\bibitem{cideal} {\sc D.A. Towers}, `C-Ideals of Lie algebras', to appear in {\em Comm. Algebra}.

\bibitem{index} {\sc D.A. Towers}, `The index complex of a maximal subalgebra of a Lie algebra'.

\end{thebibliography}
\end{document}